\documentclass[12pt]{article}
\usepackage{ifpdf}
\usepackage {graphics}
\usepackage{amsmath,mathptmx,amssymb,bm}

\usepackage{amssymb}
\usepackage{amsthm}
\newtheorem{thm}{Theorem}

\begin{document}


\title{On symmetries and conservation laws of a Gardner equation involving arbitrary functions}


\author{R. de la Rosa${}^{a}$, M.L. Gandarias${}^{b}$, M.S. Bruz\'on${}^{c}$\\
 ${}^{a}$ Universidad de C\'adiz, Spain  (e-mail: rafael.delarosa@uca.es). \\
 ${}^{b}$ Universidad de C\'adiz, Spain (e-mail:  marialuz.gandarias@uca.es). \\
 ${}^{c}$ Universidad de C\'adiz, Spain  (e-mail:  m.bruzon@uca.es). \\  
}

\date{}
 
\maketitle

\begin{abstract}

In this work we study a generalized variable-coefficient Gardner equation from the point of view of Lie symmetries in partial
differential equations. We find conservation laws by using the multipliers method of Anco and Bluman which does not require the use
of a variational principle. We also construct conservation laws using Ibragimov theorem which is based on the concept of adjoint equation for nonlinear differential equations.\\

\end{abstract}



\maketitle

\section{Introduction}
\label{Introduction}

Nonlinear equations with variable coefficients have become increasingly important in recent years because these describe many
nonlinear phenomena more realistically than equations with constant coefficients.  The Gardner equation, for instance, is used in different areas of physics, such as fluid dynamics, plasma physics, quantum field theory, and it also describes a variety of wave phenomena in plasma and solid state.\\

\noindent In this paper, we consider the
variable-coefficient Gardner equation with nonlinear terms given by

\begin{equation}\label{ed1}u_t + A(t) uu_x+ C(t) u^{2}u_x + B(t)u_{xxx} + Q(t)u =0,\end{equation}\\
where  $A(t) \neq 0$, $B(t) \neq 0$, $C(t) \neq 0$ and $Q(t)$ are arbitrary smooth functions of $t$.\\

In \cite{JoKh:10}, for $A(t)=1$ and $C(t)=0$, the optimal system of
one-dimensional subalgebras was obtained. In \cite{JoKh:11}  for
the same equation with some special forms of the functions $B(t)$
and $Q(t)$ some conservation laws were constructed. Lie symmetries
of equation (\ref{ed1}) when $Q(t)=0$, were derived in
\cite{MoRa:12}. The classification of Lie symmetries obtained in
\cite{MoRa:12} was enhanced in \cite{Sopho:14} by using the general
extended equivalence group. In \cite{ZhDoYa:08}, adding  to equation
(\ref{ed1}) the term $E(t)u_x$, where $E(t)$ is an arbitrary smooth
function of $t$, the authors found new exact non-travelling
solutions, which include soliton solutions, combined soliton
solutions, triangular periodic solutions, Jacobi elliptic function
solutions and combined Jacobi elliptic function solutions of
equation (\ref{ed1}). Soliton solutions of equation (\ref{ed1})
were obtained in \cite{wazwaz} transforming the equation to an
homogeneous equation when $Q(t)=0$ and a forcing term $R(t)$ has
been added. Finally, in \cite{HoLu:12}, exact solutions were
obtained by using the general mapping deformation method adding a
new term $E(t)u_x$ and a forcing term $R(t)$ to equation
(\ref{ed1}), where $E(t)$ and $R(t)$ are
arbitrary smooth functions of $t$.\\

Lie symmetries, in general, symmetry groups, have several applications in the context of nonlinear differential equations. It is noteworthy that they are used to obtain exact solutions and conservation laws of partial differential equations \cite{Bru:12,Bruz:14, Rafa:14, rasin,xu}.\\

 In \cite{anco} Anco and Bluman gave  a general treatment of a direct conservation law method for partial differential equations expressed
in a standard Cauchy-Kovaleskaya form
$$u_t=G(x,u,u_x,u_{xx},\ldots,u_{nx}).$$
Nontrivial conservation laws are characterized by a multiplier
$\lambda$, which has no dependence on $u_t$ and all derivatives of $u_t$, satisfying
$$\hat{E}[u]\left(\lambda u_t-\lambda G(x,u,u_x,u_{xx},\ldots,u_{nx})\right)=0.$$
Here $$\hat{E}[u]:=\frac{\partial}{\partial u}-D_t
\frac{\partial}{\partial u_t}-D_x  \frac{\partial}{\partial
u_x}+D_x^2  \frac{\partial}{\partial u_{xx}}+\ldots.$$ The
conserved density $T^t$ must satisfy
$$\lambda=\hat{E}[u] T^t,$$ and the flux $T^x$ is given by
$$ T^x=-D_x^{-1}(\lambda G)-\frac{\partial  T^t}{\partial u_x}G+GD_x\left( \frac{\partial  T^t}{\partial u_{xx}}\right)+\ldots.$$

In \cite{Ibra:07}, Ibragimov introduced a general theorem on conservation laws which does not require existence of a classical Lagrangian and is based on the concept of an adjoint equation for nonlinear equations. In \cite{Ibran:11}, Ibragimov generalized the concept of linear self-adjointness by introducing the concept of nonlinear self-adjointness of differential equations. This concept has been recently used for constructing conservation laws \cite{Rafa:14, Tra:14}.\\

The aim of this work is to obtain Lie symmetries of equation (\ref{ed1}) and construct conservation laws using both methods, the direct method proposed by Anco and Bluman \cite{anco1, anco} and Ibragimov theorem \cite{Ibra:07}. We have studied Lie symmetries of equation (\ref{ed1}) for cases $Q \neq 0$ and $Q=0$. In order to obtain conservation laws using Ibragimov theory we have determined the subclasses of equation (\ref{ed1}) which are nonlinearly self-adjoint.


\section{Classical Symmetries}
To apply the classical method to equation (\ref{ed1})  we consider
the one-parameter Lie group of infinitesimal transformations in
$(x,t,u)$ given by
$$\begin{array}{lcr}\nonumber x^* & = & x+\epsilon \xi(x,t,u)+O(\epsilon^2),
\\\nonumber t^* & = & t+\epsilon \tau(x,t,u)+O(\epsilon ^2),\\\nonumber u^* & = & u+\epsilon
\eta(x,t,u)+O(\epsilon ^2),\end{array}$$ where  $\epsilon$ is the
group parameter. We require that this transformation leaves
invariant the set of solutions of equation (\ref{ed1}). This yields to an
overdetermined, linear system of differential equations for the infinitesimals
$\xi(x,t,u),$ $\tau(x,t,u)$ and $\eta(x,t,u).$  The associated Lie
algebra of infinitesimal symmetries is formed by the set of vector fields of
the form
\begin{equation}\label{vect1}{\bf v}=\xi(x,t,u) \partial_x+
\tau(x,t,u)\partial_t+\eta(x,t,u)\partial_u.\end{equation}
Invariance of  equation (\ref{ed1}) under a Lie group of point
transformations with infinitesimal generator (\ref{vect1}) leads to
a set of 18 determining equations. By simplifying this system we obtain that $\xi=\xi(x,t)$, $\tau=\tau(t)$, and
$\eta=\eta(x,t,u)$ are related by the
following conditions:
\\
\begin{equation}\label{sis}\begin{array}{r}

\eta_{uu}=0,\\

\eta_{ux} -\xi_{xx}=0,\\

\eta_{uuu}=0,\\

\eta_{uux}=0,\\

-\tau B_t-\tau_t B+3 \xi_x B=0,\\

\tau u B Q_t-\tau u B_t Q- \eta_u u B Q+ 3 \xi_x u B Q +\eta B Q
+\eta_x u^{2} B C   \, \\   + \eta_{xxx} B^2 + \eta_x u A B+ \eta_t B=0, \\

\tau u^{3} B C_t -\tau u^{3} B_t C + \xi_x u^{3} B C
+ 2 \eta u^{2} B C- \tau u^{2} A B_t+3 \eta_{uxx} u B^2    \, \\ - \xi_{xxx} u B^2   + \tau u^{2} A_t B  + 2 \xi_x u^{2} A B +  \eta u A B-\xi_t u B=0.\\

\end{array} \end{equation}

\noindent In order to find Lie symmetries of the equation, we distinguish two cases: $Q \neq 0$ and $Q=0$.\\

\noindent{\bf Case 1.}  $Q \neq 0$.\\

\noindent For the sake of simplicity, in Case 1 we shall consider $C(t)=1$, obtaining the following symetries
$$\xi={k_1 x + \beta},\qquad \tau=\tau(t),\qquad \eta= \frac{\beta_t}{A}+\left( k_1+\alpha \right)u $$
where $A=A(t)$, $B=B(t)$, $Q=Q(t)$, $\alpha = \alpha(t)$,
$\beta =\beta (t)$ and $\tau=\tau(t)$ must satisfy the following
conditions:
\begin{eqnarray}{}
\label{eq1c2} \left(3 k_1 -\tau_{t}\right) B - \tau B_t &=&0,\\
\label{eq2c2} -  \tau B_t + \left( 4 k_1 + 2 \alpha \right)B &=&0,\\
\label{eq3c2} 2 \beta_t B  +  \tau A A_t B -  \tau A^2 B_t + \left( 3 k_1  +  \alpha \right)  A^2 B &=&0,\\
\label{eq4c2}\tau B Q_t -  \tau B_t Q + 3 k_1 B Q +  \alpha_t B  &=&0, \\ \label{eq5c2}
\beta_t \left( A Q - A_t \right) + \beta_{tt}A  &=&0.
\end{eqnarray}

\noindent Equation (\ref{eq1c2}) can be written as:

\begin{equation}\label{eqBtau} \displaystyle\frac{B_t}{B}=\frac{3k_1-\tau_t}{\tau}=k,\end{equation}

\noindent where $k$ is a constant. From (\ref{eq5c2}) we get:

\begin{equation}\label{eqDs} \displaystyle Q=\frac{A_t}{A}-\frac{\beta_{tt}}{\beta_t}. \end{equation}\\

\noindent{\bf Subcase 1.1.} Setting $k \neq 0$, solving (\ref{eqBtau}) we obtain

\begin{equation}\label{eqBts}B={ b_0}\,e^{k\,t},\qquad \tau(t)= {k_4}\,e^ {- k\,t }+{{3\,{k_1}}\over{k}}.\end{equation}\\

\noindent From (\ref{eq2c2}) we obtain


\begin{equation}\label{eqalpha}\alpha(t)= \frac{ k{k_4} e^ {- k\,t }-{k_1}  }{2}.\end{equation}\\

\noindent Substituting (\ref{eqDs}), (\ref{eqBts}) and (\ref{eqalpha}) into equations (\ref{eq3c2}) and (\ref{eq4c2}), we get the following system

\begin{eqnarray}{}
\label{eq1c11f}  - 2 A A_t f_1 + \left(k k_1 e^{k t} +k^2 k_4 \right) A^2 - 4 k e^{k t} \beta_t  &=&0,\\
\label{eq2c11f}   2 \displaystyle{\left(\frac{A_t}{A}-\frac{\beta_{tt}}{\beta_t} \right)_t} f_1 - 2 k^2 k_4 \left( \displaystyle{\frac{A_t}{A}-\frac{\beta_{tt}}{\beta_t}} \right)  -k^3 k_4&=&0.\\
\nonumber
\end{eqnarray}

\noindent where $f_1(t)=3 k_1 e^{k t}+k k_4$. This system admits two solutions

\begin{equation}\label{Apc2} A(t)= \frac{\sqrt{e^{k t}}}{2  d_0+k  k_1} f_1{}^{-\frac{d_0}{3 k k_1}-\frac{1}{2}} \left(2  a_1  k
    f_1{}^{\frac{d_0}{3 k k_1}+\frac{1}{6}}+a_0 \left(2  d_0+k
    k_1\right)\right) ,\end{equation}

\begin{equation}\label{betapc2} \beta(t)= \frac{a_0}{8 d_0 k \left(k k_1-2 d_0\right)} f_1{}^{-\frac{2 d_0}{3 k k_1}} \left(a_0 \left(4 d_0^2-k^2
   k_1^2\right)+8 a_1 d_0 k f_1{}^{\frac{d_0}{3 k
   k_1}+\frac{1}{6}}\right) ,\end{equation}

\normalsize

\noindent and

\begin{equation}\label{Anc2} A(t)= \frac{\sqrt{e^{k t}}}{2 d_0+k k_1} f_1{}^{-\frac{d_0}{3 k k_1}-\frac{1}{2}} \left(a_0
   \left(2 d_0+k k_1\right) f_1{}^{\frac{d_0}{3 k k_1}+\frac{1}{6}}-2 a_1
   k\right) ,\end{equation}

\begin{equation}\label{betanc2} \beta(t)=\frac{a_1}{2  \left(4 d_0^2- k^2 k_1^2\right)} f_1{}^{-\frac{2 d_0}{3 k k_1}} \left(2 a_0 \left(2 d_0+k
   k_1\right) f_1{}^{\frac{d_0}{3 k k_1}+\frac{1}{6}}+ \frac{a_1 k}{d_0} \left(k k_1-2
   d_0\right)\right) .\end{equation}\\

\normalsize

\noindent Lastly, substituting  solution
(\ref{Apc2})-(\ref{betapc2}) or  solution
(\ref{Anc2})-(\ref{betanc2}) into (\ref{eqDs}), we obtain

\begin{equation}\label{Dsc11} Q(t)= \frac{2 d_0 e^{k t}-k^2 k_4}{2 f_1} .\end{equation}\\

\noindent In the above equations, the following appointments are introduced: $a_0 \neq 0$, $a_1 \neq 0$, $b_0 \neq 0$, $d_0 \neq 0$, $k_1 \neq 0$, $k_4$ are arbitrary constants, $2 d_0\pm k k_1 \neq 0$. \\

%

\noindent{\bf Subcase 1.2.} Setting $k = 0$, solving (\ref{eqBtau}) we obtain

\begin{equation}\label{eqBts2}B={ b_0},\qquad \tau(t)= 3 k_1 t+k_3.\end{equation}\\

\noindent From (\ref{eq2c2}) we obtain


\begin{equation}\label{eqalpha2}\alpha(t)= -2 k_1.\end{equation}\\

\noindent Substituting (\ref{eqBts2}), (\ref{eqalpha2}) and (\ref{eqDs}) into equations (\ref{eq3c2}) and (\ref{eq4c2}), we get the following system

\normalsize

\begin{eqnarray}{}
\label{eq1c12f}   - 2 A A_t \tau -2 k_1 A^2- 4  \beta_t &=&0,\\
\label{eq2c12f}   2 \tau \displaystyle{\left(\frac{A_t}{A}-\frac{\beta_{tt}}{\beta_t} \right)_t}
 + 6 k_1  \displaystyle{\left(\frac{A_t}{A}-\frac{\beta_{tt}}{\beta_t} \right)}&=&0.\\
\nonumber
\end{eqnarray}

\normalsize

\noindent In this case, we have

\begin{equation}\label{Ac12} A(t)=a_0 \tau{}^{-\frac{d_0}{3 k_1}}+a_1 \tau{}^{-\frac{1}{3}} ,\end{equation}

\begin{equation}\label{betac12}  \beta(t)= \frac{1}{2} a_0 \left(d_0-k_1\right) \tau{}^{-\frac{2 d_0}{3 k_1}-\frac{1}{3}}
   \left(\frac{a_0 \tau{}^{\frac{4}{3}}}{3 k_1-2 d_0}-\frac{a_1 \tau{}^{\frac{d_0}{3 k_1}+1}}{d_0-2 k_1}\right) .\end{equation}\\

\normalsize \noindent Substituting  solution
(\ref{Ac12})-(\ref{betac12}) into (\ref{eqDs}), we get that

\begin{equation}\label{Ds2c12} Q(t)= \frac{d_0}{\tau} .\end{equation}\\

\noindent In the above equations, the following appointments are introduced: $a_0 \neq 0$, $a_1$, $b_0 \neq 0$, $d_0 \neq 0$, $k_1 \neq 0$, $k_3$ are arbitrary constants, $ \displaystyle{d_0\neq \frac{3}{2}k_1}$ and $d_0 \neq 2 k_1$.\\

\noindent{\bf Case 2.}  $Q=0$.\\

\noindent Now, the generators are given by:

$$\xi={k_1 x + \beta},\qquad \tau=\tau(t),\qquad \eta= \left( k_1+ k_3 \right)u +k_2,$$
where $A=A(t)$, $B=B(t)$, $C=C(t)$, $\beta =\beta (t)$ and $\tau=\tau(t)$ must satisfy the following
conditions:
\begin{eqnarray} \label{eq1Mol}\left(3 k_1 -\tau_{t}\right) B - \tau B_t &=&0,\\\label{eq2Mol}
{  \tau} B C_{t} - {  \tau} B_t C+ \left( 4 k_1 + 2 k_3  \right) B C &=&0,\\\label{eq3Mol}
 {  \tau} A_{t} B- {  \tau} A B_t + 2 k_2 B C+ \left( 3k_1 + k_3  \right)A B  &=&0,\\ \label{eq4Mol}{  k_2} A-{  \beta}_{t}&=&0.\end{eqnarray}
\\

\noindent \textbf{Subcase 2.1.} We consider $k \neq 0$. In this case
we have that $B(t)$ and $\tau(t)$ are given by (\ref{eqBts}). The
remaining functions are defined as follows

\begin{equation}\label{eqC21} C(t)= c_0 e^{k t} f_1{}^{-\frac{2k_3}{3k_1}
   -\frac{4}{3}},\end{equation}

\begin{equation}\label{eqA21} A(t)= \frac{1}{k_1+k_3} e^{k t} f_1{}^{-\frac{2 k_3}{3 k_1}
   -\frac{4}{3}} \left(2 c_0
   k_2 + a_0 \left(k_1+k_3\right) f_1{}^{\frac{k_3}{3 k_1}+ \frac{1}{3}}\right),\end{equation}

\begin{equation}\label{eqbeta21} \beta(t)=\beta
   _0 -\frac{a_0 k_2 f_1{}^{-\frac{k_3}{3 k_1}}}{k k_3}-\frac{2 c_0 k_2^2 f_1{}^{-\frac{2 k_3}{3
   k_1}-\frac{1}{3}}}{k \left(k_1+k_3\right) \left(k_1+2 k_3\right)}.\end{equation}\\

\noindent where $f_1$ has already been previously defined as $f_1(t)=3 k_1 e^{k t}+k k_4$.\\

\noindent \textbf{Subcase 2.2.} We consider $k = 0$. Now, $B(t)$ and
$\tau(t)$ are given by (\ref{eqBts2}). In this case, the remaining
functions are given by

\begin{equation}\label{eqC22} C(t)= c_0 \tau{}^{-\frac{2 k_3}{3 k_1}-\frac{4}{3}},\end{equation}

\begin{equation}\label{eqA22} A(t)=\frac{1}{k_1+k_3}\tau{}^{-\frac{2 k_3}{3 k_1}
   -\frac{4}{3}}
   \left(2 c_0 k_2 +a_0 \left(k_1+k_3\right) \tau{}^{\frac{k_3}{3 k_1}+ \frac{1}{3}}\right) ,\end{equation}

\begin{equation}\label{eqbeta22} \beta(t)= \beta
   _0 -\frac{a_0 k_2 \tau{}^{-\frac{k_3}{3 k_1}}}{k_3}-\frac{2 c_0 k_2^2 \tau{}^{-\frac{2 k_3}{3
   k_1}-\frac{1}{3}}}{\left(k_1+k_3\right) \left(k_1+2 k_3\right)}.\end{equation}\\

\noindent In Case 2,  the following designations are introduced: $a_0$, $b_0 \neq 0$, $c_0 \neq 0$, $\beta_0$, $k_1 \neq 0$, $k_2$, $k_3 \neq 0$, $k_4$ are arbitrary constants, $ k_1\neq -k_3$ and $k_1 \neq -2 k_3$. For this case the symmetries were obtained in \cite{MoRa:12}.


\section{Formal Lagrangian and adjoint equation}

In \cite{Ibra:07} Ibragimov introduced a new theorem on conservation
laws. This theorem is valid for any system of differential equations
wherein the number of equations is equal to the number of dependent
variables. The new theorem does not require existence of a classical
Lagrangian and it is based on the concept of adjoint equation
for nonlinear equations. In order to obtain the adjoint equation we use the following
definition:\\

\noindent \textbf{Definition 1.} Consider an qth-order partial differential equation
 \begin{equation}
 \label{fa}
 F({\rm x},u,u_{(1)}, \ldots,u_{(s)}) =0
 \end{equation}
 with independent variables ${\rm x}=(x^1,\ldots,x^n)$ and
 a dependent variable $u,$ where
$u_{(1)}=\{u_i\},$ $u_{(2)}=\{u_{ij}\},\ldots\,$
 denote the sets of
partial derivatives  of first, second, etc. order,
$u_i=\partial u/\partial x^i$, $u_{ij}=\partial^2 u/\partial
x^i\partial x^j.$ The formal Lagrangian is defined as
\begin{equation}\label{lag}{\cal L}=v\,F\left({\rm x},u, u_{(1)}, \ldots,u_{(s)}\right), \end{equation} where
$v=v({\rm x})$ is a new dependent variable.  The adjoint equation to
{\rm(\ref{fa})} is
 \begin{equation}
 \label{faadj}
 F^{*}({\rm x},u,v,u_{(1)}, v_{(1)},\ldots,u_{(s)},v_{(s)})
 =0, \end{equation}
 with

 $$F^{*}({\rm x},u,v,u_{(1)}, v_{(1)},\ldots,u_{(s)},v_{(s)})
 =\frac{\delta(v\,F)}{\delta u},$$
where

$$ \frac{\delta }{\delta u}=\frac{\partial}{\partial u}
 +\displaystyle\sum_{s=1}^{\infty} (-1)^s D_{i_1}\cdots
D_{i_s}\frac{\partial}{\partial u_{i_1\cdots i_s}},$$ denotes the
variational derivatives (the Euler-Lagrange
 operator).
 Here
  $$
 D_i=\displaystyle\frac{\partial}{\partial
x^i}+u_i\frac{\partial}{\partial u}+ u_{ij}\frac{\partial}{\partial
u_j}+\cdots
 $$
represents the total differentiation.

\begin{thm} The adjoint
equation to equation (\ref{ed1})  is
\begin{equation}\label{eqad}F^*\equiv  v\,Q-u^2\,v_{x}\,C-v_{xxx}\,B-u\,v_{x}\,A-v_{t}.\end{equation}
\end{thm}

\section{Nonlinearly self-adjoint equations}
\noindent In this section we use the following definition given in \cite{Ibran:11}.\\

\noindent \textbf{Definition 2.} Equation {\rm (\ref{fa})} is said to be {\bf nonlinearly self-adjoint} if the
 equation obtained from the adjoint equation {\rm (\ref{faadj})}
 by the substitution \begin{equation}\label{newv}v=\varphi({\rm x},u),\end{equation}
 such that $\varphi({\rm x},u)\neq 0$,
 is identical with the original equation {\rm(\ref{fa})}, i.e.
 \begin{equation}\label{condadj}
 F^*\,{|}_{v=\varphi}= \lambda ({\rm x}, u, \dots  ) F.\end{equation}
 for some differential function $\lambda= \lambda ({\rm x}, u, \dots  )$.
 If $\varphi=u$  or $\varphi=\varphi(u)$ and $\varphi'(u)\neq 0$, equation {\rm (\ref{fa})}
 is said {\bf self-adjoint} or {\bf quasi self-adjoint},
 respectively. If $\varphi_{x^i}({\rm x},u)\neq 0$ or $\varphi_u({\rm x},u)\neq 0$ equation {\rm (\ref{fa})}
 is said  {\bf  weak
  self-adjoint} \cite{mlfc}.\\

\noindent Taking into account expression (\ref{eqad}) and using
(\ref{newv}) and its derivatives, equation (\ref{condadj}) can be written
as

\small
\begin{equation}\label{eqlambda}\begin{array}{l}  u_{x}\,\left( \lambda \left(-u^2\,C-u\,A\right)-\varphi_{u}\,u^2\,C
 -3\,\varphi_{uu}\,u_{xx}\,B-3\,\varphi_{uxx}\,B-\varphi_{u}
 \,u\,A\right) \\ \\

 +u_{xxx}\,\left(-\lambda \, B-\varphi_{u}\,B\right)-\lambda \, u
 \,Q +u_{t}\,\left(-\lambda-\varphi_{u}\right) +\varphi\,Q-
 \varphi_{x}\,u^2\,C -\varphi_{t}  \\ \\

 -3\,\varphi_{ux}\,u_{xx}\,B-\varphi_{uuu}
 \,\left(u_{x}\right)^3\,B-3\,\varphi_{uux}\,\left(u_{x}\right)^2
 \,B-\varphi_{xxx}\,B-\varphi_{x}\,u\,A =0  .\end{array}\end{equation}

\normalsize

\noindent Equation (\ref{eqlambda}) should be satisfied identically in all
variables $u_t$, $u_{x}$, $u_{xx},\ldots$ Equating to zero the
coefficients of the derivatives of $u$ we obtain:

\begin{thm} Equation (\ref{ed1}) with $A(t) \neq 0, B(t)\neq 0$, $C(t) \neq 0$
and $Q(t)$ arbitrary functions, is nonlinearly self-adjoint and
\begin{equation}\label{hs}\varphi= \displaystyle c_1 e^{\int 2 Q(t) dt} u+ c_2 e^{\int  Q(t) dt},\end{equation}
with $c_1$ and $c_2$ arbitrary constants.
\end{thm}

\section{Conservation Laws}

Conservation laws appear in many of physical, chemical and mechanical processes, such laws enable us to solve problems in which certain physical properties do not change over time within an isolated physical system.\\

\noindent A conservation law of equation (\ref{ed1}) is a space-time divergence such that
\medskip
\begin{equation}\label{cl}
 \begin{array}{rcc}
\displaystyle{D_t T^t(x,t,u,u_x,u_t,...)+D_x T^x(x,t,u,u_x,u_t,...)=0,}\\
\end{array}
\\
\end{equation}\\
on all solutions $u(x,t)$ of the equation (\ref{ed1}). Here, $T^t$ represents the conserved density and $T^x$ the associated flux \cite{ancoibra}, and $D_x$, $D_t$ denote the total derivative operators with respect to $x$ and $t$ respectively.\\

\noindent In this section we construct conservation laws of each case by using both methods.

\subsection{Conservation laws by using the direct method of the multipliers of Anco and Bluman}

\noindent We suppose,  that $T^t$ and $T^x$ have no dependence on $u_t$ and all derivatives of $u_t$.\\
\\
Each conservation law (\ref{cl}) has an equivalent characteristic form in which has been eliminated $u_t$ and its differential consequences from $T^t$ and $T^x$ by using equation (\ref{ed1})
\begin{eqnarray*}
 \widehat{T^t}= T^t \mid_{u_t=\Delta}= T^t- \Phi,  \\ \\
 \widehat{T^x}= T^x \mid_{u_t=\Delta}= T^x- \Psi,
\end{eqnarray*}
\noindent where $ \Delta= - A u u_x- C u^2 u_x-B u_{xxx} -Q u$, so that
$$ \left(D_t \widehat{T^t}(x,t,u,u_x,u_{xx},...)+D_x \widehat{T^x}(x,t,u,u_x,u_{xx},...)\right)\mid_{u_t=\Delta}=0, $$
\noindent is verified on all solutions of equation (\ref{ed1}), and where
\begin{eqnarray*}{} D_t \mid_{u_t=\Delta} &=& \partial_t+ \Delta\partial_u +D_{x}(\Delta)\partial_{u_x}+...    \\ \\
 D_x \mid_{u_t=\Delta} & = &  \partial_x+ u_x \partial_u+ u_{xx} \partial_{u_x}+...=D_x.
 \end{eqnarray*}
\noindent In particular, moving off of solutions, we have the identity
\begin{eqnarray*} D_t = D_t \mid_{u_t=\Delta} + (u_t + Auu_x+ Cu^{2}u_x + Bu_{xxx} + Qu)\partial_u \quad \quad \\ \\  + D_x(u_t + Auu_x+ Cu^{2}u_x + Bu_{xxx} + Qu)\partial_{u_x}+... \end{eqnarray*}
These expressions yield us to the characteristic form of
conservation law (\ref{cl}) \small
\begin{equation}\label{char} \begin{array}{r}  D_t \widehat{T^t}(x,t,u,u_x,u_{xx},...)+D_x \left(\widehat{T^x}(x,t,u,u_x,u_{xx},...)+\widehat{\Psi}(x,t,u,u_x,u_t,...) \right)  \\ \\ = \left( u_t +A uu_x+ C u^{2}u_x + B u_{xxx} + Q u \right) \Lambda(x,t,u,u_x,u_{xx},...), \end{array}\end{equation}
\normalsize
where
\begin{eqnarray*} \widehat{\Psi}(x,t,u,u_x,u_t,...)= E_{u_x}(\widehat{T^t})\left( u_t + Auu_x+ Cu^{2}u_x + Bu_{xxx} + Q u \right) \\ \\ +E_{u_{xx}}(\widehat{T^t})D_x\left( u_t + Auu_x+ C u^{2}u_x + Bu_{xxx} + Q u \right)+... \end{eqnarray*}
is a trivial flux \cite{ancoibra}, and the function
$$ \Lambda=E_u (\widehat{T^t}),$$
is a multiplier, where $E_u= \partial_u-D_x \partial_{u_x}+D_x^{2}\partial_{u_{xx}}-...$, denotes the (spatial) Euler operator with respect to u.
\\ \\
A function $\Lambda(x,t,u,u_x,u_{xx},...)$ is called multiplier if it verifies that $(u_t+Au u_x+ Cu^2 u_x+Bu_{xxx}+ Qu)\Lambda$ is a divergence expression for all functions $u(x,t)$, not only solutions of equation (\ref{ed1}).
\\ \\
In order to obtain conservation laws, we use (\ref{char}) from which
is deduced that all nontrivial conserved densities in the form
(\ref{cl}) arise from multipliers $\Lambda$ of equation
(\ref{ed1}), where $\Lambda$ depends only on $x$, $t$, $u$ and $x$
derivatives of $u$. Divergence condition can be characterized as
follows

\begin{equation}\label{var} \frac{\delta}{\delta u} \left(  \left(   u_t + Auu_x+ Cu^{2}u_x + Bu_{xxx} + Qu  \right) \Lambda  \right) =0, \end{equation}\\

\noindent where $ \displaystyle{ \frac{\delta}{\delta u}= \partial_u-D_x \partial_{u_x}-D_t \partial_{u_t}+D_x D_t\partial_{u_{xt}}+D_x^2\partial_{u_{xx}}+...} $, denotes the variational derivative.
\\
\\
Equation (\ref{var}) is linear in $u_t$, $u_{tx}$, $u_{txx}$,... taking the coefficient of $u_t$ and $x$ derivatives of $u_t$ we obtain a system of determining equations for $\Lambda $, which yields the following equivalent equations \cite{anco0,anco1,anco}
$$ -D_t \Lambda-  D_x \left( \left( Au+Cu^2 \right) \Lambda \right)- D_x^3 \left(B \Lambda \right)+ \left(A u_x +2 C u u_x+ Q \right) \Lambda=0,$$
and
\\
$$ \Lambda_u= E_u(\Lambda), \qquad \, \Lambda_{u_x}=-E_u^{(1)}(\Lambda), \qquad \, \Lambda_{u_{xx}}= -E_u^{(2)}(\Lambda), ... $$\\
which are verified for all solutions $u(x,t)$ of equation (\ref{ed1}).
\\

\noindent Given a multiplier $\Lambda$, we can obtain the conserved density using a standard method \cite{wolf}
$$ T^t= \int_{0}^{1} d \lambda \, \, u \Lambda (x,t,\lambda u,\lambda u_x,\lambda u_{xx},...).$$
\\
We have considered multipliers up to second order, i.e., $\Lambda(x,t,u,u_x,u_{xx})$. In this section we proceed to obtain conservation laws for the values obtained for $A(t)$, $B(t)$, $C(t)$ and $Q(t)$ in Case 1, Section 2.\\

\noindent \textbf{Subcase 1.1.} In this case, equation (\ref{ed1}) is given by

\begin{equation}\label{subcase11}\displaystyle u_{{t}}+A(t) uu_{{x}} +{u}^{2}u_{{x}}+{ b_0}{{ e}^{kt}}u_{xxx}+{\frac{\ 2{ d_0}{{ e}^{kt}}-{k}^{2}{
 k_4}  }{2 f_1}}u
 =0 , \end{equation}

\noindent where $A(t)$ is given by (\ref{Apc2}) or (\ref{Anc2}). For
equation (\ref{subcase11})  multiplier is given by

$$ \Lambda= \displaystyle e^{-\frac{k t}{2}} f_1{}^{{\frac {2 d_{{0}}}{3 k k_{{1}}}+1}
}\left(  c_1 +c_2 \, e^{-\frac{k t}{2} } \, u \right),$$

\noindent where $f_1(t)=3 k_1 e^{k t}+k k_4$. The conserved density obtained from this multiplier is

$$ \begin{array}{lll} T^t & = & \displaystyle e^{-\frac{k t}{2}} f_1{}^{{\frac {2 d_{{0}}}{3 k k_{{1}}}+1}
} \, u \left(  c_1 +\frac{c_2}{2} \, e^{-\frac{k t}{2} } \, u \right). \end{array}
$$

\noindent The flux obtained from this multiplier depends on function
$A(t)$. If $A(t)$ is given by (\ref{Apc2}), we obtain

\normalsize

$$ \begin{array}{lll} T^x & = & \displaystyle \frac{{\rm e}^{-\frac{kt}{2}}}{12} f_1^{\frac{d_0}{3 k k_1}}
 \left( 3 \, c_2 \, f_1{}^{\frac{ d_0}{3 k k_1}+1} \left(  u^4 \, e^{-\frac{kt}{2}} +2 \, b_0  \, e^{\frac{k t}{2}} \left( 2 \, u \, u_{xx}- u_x^2 \right)     \right)
\right. \\ \\

 & & +4 \left( a_0 \, c_2+c_1 \right) u^3 \,  f_1 + 6 \, a_0 \, c_1 \, u^2 \, e^{\frac{k t}{2}}f_1^{-\frac{d_0}{3 k k_1}}

 +12 \, b_0 \, c_1 \, u_{xx} \, e^{k t} \, f_1^{\frac{1}{2}}\\ \\

& & +\displaystyle \left. \frac{4 \, a_1 \, k }{2 \, d_0 + k \, k_1} u^2 \left( 2 \, c_2  \, u \, f_1^{\frac{d_0}{3 k k_1}+\frac{2}{3}} + 3  \, c_1 \, e^{\frac{k t}{2}} \, f_1^{\frac{1}{6}}  \right)
 \right) .
 \end{array}$$

\normalsize

\noindent If $A(t)$ is given by (\ref{Anc2}), we get

$$ \begin{array}{lll} T^x & = & \displaystyle \frac{{\rm e}^{-\frac{kt}{2}}}{12} f_1^{\frac{d_0}{3 k k_1}}
 \left( 3 \, c_2 \, f_1{}^{\frac{ d_0}{3 k k_1}+1} \left(  u^4 \, e^{-\frac{kt}{2}} +2 \, b_0  \, e^{\frac{k t}{2}} \left( 2 \, u \, u_{xx}- u_x^2 \right)     \right)
\right. \\ \\

 & & +4 \, c_1 \, f_1^{\frac{1}{2}} \left(  u^3 + 3 \, b_0 \, u_{xx} \, e^{k t}   \right) + 6 \, a_0 \, c_1 \, u^2 \, e^{\frac{k t}{2}}f_1^{\frac{1}{6}}

+ 4 \, a_0 \, c_2 \, u^3 \, f_1^{\frac{d_0}{3 k k_1}+ \frac{2}{3}}\\ \\

& & - \displaystyle \left. \frac{4 \, a_1 \, k }{2 \, d_0 + k \, k_1} u^2 \left( 2 \, c_2  \, u \, f_1^{\frac{1}{2}} + 3  \, c_1 \, e^{\frac{k t}{2}}  \right)
 \right) .
 \end{array}$$\\

\normalsize

\noindent \textbf{Subcase 1.2.} In this case, equation (\ref{ed1}) is given by

\normalsize
\begin{equation}\label{subcase12}    \displaystyle u_t + \left( a_0 \tau{}^{-\frac{d_0}{3 k_1}}+a_1 \tau{}^{-\frac{1}{3}} \right) u u_x+
 u ^{2}u_x +b_{{0}}u_{xxx} +\frac{d_0}{\tau}u=0, \end{equation}

\noindent where $\tau$ is given by (\ref{eqBts2}). For equation (\ref{subcase12}) has been obtained the following multiplier

$$ \Lambda=  \tau{}^{{\frac {d_{{0}}}{3 k_{{1}}}}
} \left( c_1 +c_2 \, \tau{}^{{\frac {d_{{0}}}{3 k_{{1}}}}
}\right).$$\\

\noindent The conserved density and the flux obtained from this multiplier are:
\\
$$ \begin{array}{lll} T^t & = & \displaystyle \tau{}^{{\frac {d_{{0}}}{3 k_{{1}}}}
} \, u \left( c_1 + \frac{c_2}{2}  \, \tau{}^{{\frac {d_{{0}}}{3 k_{{1}}}}
} \, u \right) , \end{array}
$$\\
$$ \begin{array}{lll} T^x & = & \displaystyle   {\frac {1}{12}} \, \tau{}^{{\frac {d_{{0}}}{3 k_{{1}}}}
} \left( \tau{}^{{\frac {d_{{0}}}{3 k_{{1}}}}
} \left(   6 \, b_0 \, c_2 \left( 2 \, u \, u_{xx} -u_x^2 \right) +3 \, c_2 \, u^4  \right) +6 \, a_0 \, c_1 \, u^2 \, \tau^{-\frac{d_0}{3 k_1}} \right. \\ \\

 & & \left. +6 \, a_1 \, c_1 \, u^2 \, \tau^{-\frac{1}{3}}+4 \, a_1 \, c_2 \, u^3 \, \tau^{\frac{d_0}{3 k_1}-\frac{1}{3}} + 4 \, c_1 \left(      u^3 + 3 \, b_0 \, u_{xx}\right) +4 \, a_0 \, c_2 \, u^3 \right)
.  \end{array}$$\\

\subsection{Conservation laws by using a general theorem on conservation laws proved by Ibragimov}

In this section we construct conservation laws for the values obtained for the functions $A(t)$, $B(t)$, $C(t)$ and $Q(t)$ in Case 2, Section 2 using the following theorem on conservation laws proved in \cite{Ibra:07}.

\begin{thm} Any Lie point, Lie-B\"{a}cklund or non-local symmetry
 \begin{equation}
 \label{gee}
 {\bf v}=\xi^i({\rm x},u,u_{(1)},\ldots)\frac{\partial}{\partial x^i}
 +\eta({\rm x},u,u_{(1)},\ldots)\frac{\partial}{\partial u},
 \end{equation} of  equation {\rm (\ref{fa})}
 provides a conservation law $\,D_i(T^i)=0\, $ for the simultaneous system {\rm (\ref{fa})}, {\rm
 (\ref{faadj})}. The conserved vector is given by
  \begin{equation}
 \label{e19}
 \begin{array}{ll}
 T^i&=\xi^i{\cal L} +
 W\left[\displaystyle\frac{\partial{\cal L}}{\partial
 u_i}-D_j\left(\frac{\partial {\cal L}}{\partial
 u_{ij}}\right)+D_jD_k\left(\frac{\partial {\cal
 L}}{\partial u_{ijk}}\right)-\cdots\right]\\[3.5ex] &
 +D_j(W)\left[\displaystyle\frac{\partial {\cal
 L}}{\partial u_{ij}}-D_k\left(\frac{\partial {\cal L}}{\partial
 u_{ijk}}\right)+\cdots\right]
 +D_jD_k(W)\left[\displaystyle\frac{\partial
 {\cal L}}{\partial u_{ijk}}-\cdots\right]+\cdots\,,
 \end{array}
 \end{equation}
 where $W$ and ${\cal L}$ are defined as follows:
 \begin{equation}
 \label{la}
 W = \eta - \xi^j u_j, \quad
 {\cal L}=v\,F\left({\rm x},u, u_{(1)}, \ldots,u_{(s)}\right).
 \end{equation}
 \end{thm}

\noindent In order to apply Theorem 3 to our equation we perform the following change of notation:
$$\left(  x^1, x^2  \right)= \left(  t, x    \right), \qquad \, \, \left(  \xi^1, \xi^2  \right)= \left(  \tau , \xi    \right), \qquad \, \,  \left(  T^1, T^2  \right)= \left(  T^t, T^x    \right).$$ \\

\noindent \textbf{Subcase 2.1.} In this case, we consider equation (\ref{ed1}) with $Q(t)=0$, $A(t)$, $B(t)$ and $C(t)$ given by (\ref{eqA21}) , (\ref{eqBts}) and (\ref{eqC21}) respectively. Consequently, the equation admits the following generator

\begin{equation}\label{generator21}
v= \left(k_1 x + \beta(t)\right) \partial_x+ \left( k_4 e^{-k t}+ \frac{3 k_1}{k} \right) \partial_t + \left( (k_1+k_3)u+k_2 \right)\partial_u,
\end{equation}

\noindent where $\beta(t)$ is given by (\ref{eqbeta21}). From
Theorem 2, we have that $\varphi$ is given by
\begin{equation}
\label{phi21} \varphi= c_1 u+ c_2.
\end{equation}
\noindent Thus, we obtain conservation law (\ref{cl}) with the conserved vector

\begin{equation}\label{lc1} \begin{array}{lll} T^t & = & \displaystyle \left( k_3 \, u+k_2  \right) \left( c_1 \, u+c_2  \right)+ k_1 \, u \left(  \frac{3}{2}\,c_1 \, u+2 \, c_2 \right) , \end{array}
\end{equation}


$$ \begin{array}{lll} T^x & = & \displaystyle \frac{e^{k\,t}}{12\,\left(k_{3}+k_{1}\right)}   \left(     6 \, b_0 \, c_1 \left( 2 \, u \, u_{xx} -u_x^2 \right) \left( 2 \,k_3^2+5 \, k_1 \, k_3 + 3 \, k_1^2 \right)  \right. \\ \\

 & & \displaystyle +12 \, b_0 \, u_{xx} \left( c_2 \left( k_3^2 + 3\, k_1 \, k_3+ 2 \, k_1^2 \right) +c_1 \left( k_3 +k_1 \right) \right) \\ \\

 & & \displaystyle + f_1^{-\frac{k_3}{3 \, k_1}-1} \left(  4 \, a_0 \, c_1 \, u^3 \left( 2 \,k_3^2+5 \, k_1 \, k_3 + 3 \, k_1^2 \right) \right.  \\ \\

 & & \left. +6 \, a_0 \, c_2 \, u^2 \left( k_3^2+3 \, k_1 \, k_3 + 2 \, k_1^2  \right) + 6 \, a_0 \, c_1 \, k_2 \, u^2 \left( k_3 + k_1 \right)   \right) \\ \\

 & &  + f_1^{-\frac{2 \, k_3}{3 \, k_1}-\frac{4}{3}} \left( 3\, c_0 \, c_1 \, u^4 \left( 2 \,k_3^2+5 \, k_1 \, k_3 + 3 \, k_1^2 \right) + 4 \, c_0 \, c_1 \, k_1 \, u^3 \left( 5 \, k_3 + 7\, k_2 \right)  \right. \\ \\

 & &  \left. \left. +4 \, c_0 \, c_2 \, u^3 \left( k_3^2 + 3\, k_1 \, k_3+ 2 \, k_1^2 \right) + 12 \, c_0 \, k_2 \, u^2
 \left(  c_2 \, k_3 + c_1 \, k_2 + 2 \, c_2 \, k_1 \right) \right) \right).

\end{array}$$\\

\noindent \textbf{Subcase 2.2.} Now, we consider equation (\ref{ed1}) with $Q(t)=0$, $A(t)$, $B(t)$ and $C(t)$ given by (\ref{eqA22}) , (\ref{eqBts2}) and (\ref{eqC22}) respectively. Consequently, the equation admits the following generator

\begin{equation}\label{generator22}
v= \left(k_1 x + \beta(t)\right) \partial_x+ \tau \partial_t + \left( (k_1+k_3)u+k_2 \right)\partial_u,
\end{equation}

\noindent where $\tau(t)$ and $\beta(t)$ are given by (\ref{eqBts2}) and (\ref{eqbeta22}) respectively. As in the previous case, $\varphi$ is given by (\ref{phi21}). Thus, we obtain conservation law (\ref{cl}) with the conserved vector

$$ \begin{array}{lll} T^t & = & \displaystyle \left( k_3 \, u+k_2  \right) \left( c_1 \, u+c_2  \right)+ k_1 \, u \left(  \frac{3}{2}\,c_1 \, u+2 \, c_2 \right) , \end{array}
$$

$$ \begin{array}{lll} T^x & = & \displaystyle \frac{1}{12\,\left(k_{3}+k_{1}\right)}   \left(     6 \, b_0 \, c_1 \left( 2 \, u \, u_{xx} -u_x^2 \right) \left( 2 \,k_3^2+5 \, k_1 \, k_3 + 3 \, k_1^2 \right)  \right. \\ \\

 & & \displaystyle +12 \, b_0 \, u_{xx} \left( c_2 \left( k_3^2 + 3\, k_1 \, k_3+ 2 \, k_1^2 \right) +c_1 \left( k_3 +k_1 \right) \right) \\ \\

 & & \displaystyle + \tau^{-\frac{k_3}{3 \, k_1}-1} \left(  4 \, a_0 \, c_1 \, u^3 \left( 2 \,k_3^2+5 \, k_1 \, k_3 + 3 \, k_1^2 \right) \right.  \\ \\

 & & \left. +6 \, a_0 \, c_2 \, u^2 \left( k_3^2+3 \, k_1 \, k_3 + 2 \, k_1^2  \right) + 6 \, a_0 \, c_1 \, k_2 \, u^2 \left( k_3 + k_1 \right)   \right) \\ \\

 & &  + \tau^{-\frac{2 \, k_3}{3 \, k_1}-\frac{4}{3}} \left( 3\, c_0 \, c_1 \, u^4 \left( 2 \,k_3^2+5 \, k_1 \, k_3 + 3 \, k_1^2 \right) + 4 \, c_0 \, c_1 \, k_1 \, u^3 \left( 5 \, k_3 + 7\, k_2 \right)  \right. \\ \\

 & &  \left. \left. +4 \, c_0 \, c_2 \, u^3 \left( k_3^2 + 3\, k_1 \, k_3+ 2 \, k_1^2 \right) + 12 \, c_0 \, k_2 \, u^2
 \left(  c_2 \, k_3 + c_1 \, k_2 + 2 \, c_2 \, k_1 \right) \right) \right).

\end{array}$$\\

\noindent We remark that some of these conservation laws yield conserved
integrals with physical meaning. Setting in (\ref{lc1})
$k_1=1$, $c_2=\displaystyle \frac{1}{2}$ and $k_2=k_3=c_1=0$ we have
$$C_1=\int_{-\infty}^{\infty}{u \, dx},$$ which is the conserved mass for
equation (\ref{ed1}). The conserved integral arising from
(\ref{lc1}) setting $k_1=1$, $c_1=\displaystyle \frac{2}{3}$ and $k_2=k_3=c_2=0$ gives the energy
$$C_2=\int_{-\infty}^{\infty}{u^2dx},$$ for
equation (\ref{ed1}).\\

\noindent Another powerful application of conservation laws taking into
account the relationship between Lie symmetries and conservation
laws it is the so called double reduction method given by Sj\"{o}berg \cite{sjoberg}. This
method allow us to reduce the Gardner equation to a second order
ordinary differential equation. Sj\"{o}berg
introduced this method in order to get solutions of a {\textit{q}}th
partial differential equation (\ref{fa}) from the solutions of an ordinary
differential equation of order \textit{q}$-1$. This method can be applied when a symmetry {\bf v} is associated to a conserved vector $T$. In accordance with the definition given by Kara and Mahomed \cite{kara} we will establish that a symmetry {\bf v} is associated to $T$ if the following equation holds

\begin{equation}\label{associated}
{\bf v} \left( T^i \right) + T^i D_k \left( \xi^k \right) - T^k D_k \left( \xi^i \right)=0.
\end{equation}

\noindent In the terms of the canonical variables $r$, $s$ and $w$, symmetry (\ref{vect1}) becomes a translation on s, ${\bf v}=\displaystyle \frac{\partial}{\partial s}$. Thus, the conservation law can be rewritten as

\begin{equation}\label{newcl} D_s T^s+ D_r T^r=0, \end{equation}

\noindent with

\begin{equation}\label{newdensity} T^s=\displaystyle \frac{T^t D_t(s)+T^x D_x(s)}{D_t(r) D_x(s)-D_x(r)D_t(s)},\end{equation}

\noindent and

\begin{equation}\label{newflux} T^r=\displaystyle \frac{T^t D_t(r)+T^x D_x(r)}{D_t(r) D_x(s)-D_x(r)D_t(s)}.\end{equation}

\noindent Due to solutions of the equation (\ref{fa}) written in canonical variables must be invariant with respect to ${\bf v}$ and $T$ is associated with ${\bf v}$, equation (\ref{newcl}) becomes\\

$D_r T^r=0,$\\

\noindent so that

\begin{equation}\label{ode}
T^r\left(r,w,w_r,w_{rr},\dots,w_{r^{q-1}}\right)=k, \qquad k=const.
\end{equation}

\noindent We stress that equation (\ref{ode}) is an ordinary differential equation of order \textit{q}$-1$, whose solutions are solutions of equation (\ref{fa}), by writting this solution in terms of $x$, $t$ and $u$.\\

\noindent In order to show the above explained procedure, let us consider the equation

\begin{equation}\label{example}
F \equiv u_t +  uu_x+  u^{2}u_x +  u_{xxx}=0,
\end{equation}

\noindent It can be easily chequed that equation
(\ref{example}) admits the following generator

\begin{equation}\label{symmetry}
{\bf v}= \displaystyle \left( k_1 \, x-\frac{k_1 }{2 } \, t \right) \partial_x+ \left( 3 \, k_1 \, t +k_2  \right) \partial_t+ \left( -k_1 \, u-\frac{k_1 }{2} \right)\partial_u,
\end{equation}

\noindent where $k_1$ and $k_2$ are arbitrary constants. From
Theorem (\ref{gee}), one can get the following conserved vector for
generator (\ref{symmetry})

\begin{equation}\label{conservedvector}
\begin{array}{lll} T^t & = & \displaystyle -\frac{c_1 \, k_1 }{2}\, u^2 - \frac{ c_1 \, k_1}{2 }\, u-\frac{ c_2 \, k_1}{2 } ,\\ \\

T^x & = & \displaystyle - c_1 \, k_1 \, u \, u_{xx}- \frac{c_1 \, k_1 }{2}\, u_{xx}+ \frac{c_1 \, k_1}{2} \, u_{x}^2 \\ \\ & & \displaystyle  - \frac{c_1 \, k_1 }{4} \, u^4 - \frac{c_1 \, k_1}{2}\, u^3 -\frac{c_1 \, k_1 }{4}\, u^2, \end{array}
\end{equation}

\noindent where $c_1$ and $c_2$ are arbitrary constants. The conserved vector (\ref{conservedvector}) is not associated to symmetry
(\ref{symmetry}). 
%
%
It can be easily seen that generator
$$ {\bf v}= c \, \partial_x+ \partial_t.$$
is associated to (\ref{conservedvector}). The canonical coordinates
are $$ r= x-c \, t, \quad s=t, \quad
w=u.
$$

\noindent We suppose without loss of generality that $k_1=c_1=1$ and $c_2=0$. From (\ref{newflux}) we get

\begin{equation} \begin{array}{lll} T^r & = & \displaystyle\frac{1}{4} \left( \left(4\,w+2\right)\,w_{rr}-2\,w_{r}^2 \right. \left.  +w^4+2\,
 w^3+\left(1-2\,c\right)\,w^
 2   -2\,c\,w \right)$$. \end{array}
\end{equation}

\noindent Setting $T^r=\displaystyle \frac{k}{4}$, $k=\mbox{const.}$, we obtain

\begin{equation}\label{odeexample} \begin{array}{l}  \displaystyle \left(4\,w+2\right)\,w_{rr}-2\,w_{r}^2 -2\,c\,w  +\left(1-2\,c\right)\,w^
 2  +2\,
 w^3 +w^4 = k$$. \end{array}
\end{equation}

\noindent Due to the fact that (\ref{odeexample}) is an autonomous
equation the substitution  $w_r=p(w)$ yields the following first
order ordinary differential equation

$$ p' = \frac{k+2\, p^2+ 2\,c\,w+ \left(2\,c- 1\right)\,w^
 2 -2\,
 w^3 -w^4}{4\,w+2},    $$
\noindent whose solutions are solutions of (\ref {example}) when written in terms of $x$, $t$ and $u$.

\section{Conclusions}

In this paper we have considered a generalized variable-coefficient Gardner equation. Classical symmetries of equation
(\ref{ed1}) have been obtained involving different arbitrary functions which can be used to determine similarity and exact solutions. Symmetries obtained in this work generalise those already obtained by other authors in equations which belong to the family of equations (\ref{ed1}), such as KdV equation and other Gardner equations with time-dependent coefficients. In particular, we have considered two different cases, which in turn would lead to different subcases.\\

We have determined the subclasses of equation (\ref{ed1}) which are
nonlinearly self-adjoint, as well as the multipliers, of Anco and
Bluman method. We have derived conservation laws by using both
methods. We have shown that some of these conservation laws yields
conserved integrals with physical meaning, such as mass and energy.
Finally, as an example of another application of the conserved
vectors, we have applied the double reduction method to get exact
solutions of the Gardner equation from solutions of a second order
reduced ordinary differential equation.

\section*{Acknowledgements}

The authors acknowledge the financial support from Junta de
Andaluc\'{i}a group FQM-201. The second and third author also
acknowledge the support of DGICYT project MTM2009-11875 and  FEDER
funds.






\end{document}